\documentclass[12pt]{amsart}
 \usepackage{amsmath} 
 \usepackage{amssymb} 
 \usepackage{latexsym} 
\newcommand{\commsqu}[8]{\begin{center} 
\begin{picture}(6,6) 
\put(1,5){$#1$} 
\put(5,5){$#2$} 
\put(5,1){$#3$} 
\put(1,1){$#4$} 
\put(2,5.2){\vector(1,0){2.5}} 
\put(2,1.3){\vector(1,0){2.5}} 
\put(1.5,4.5){\vector(0,-1){2.5}} 
\put(5.3,4.5){\vector(0,-1){2.5}} 
\put(3,5.7){$#5$} 
\put(5.7,3){$#6$} 
\put(3,1.8){$#7$} 
\put(0.5,3){$#8$} 
\end{picture} 
\end{center}} 
\newcommand{\brref}[1]{(\ref{#1})} 
\newcommand{\tensor}{\otimes}

\newcommand{\restrict}[2]{{#1}_{\mid _{#2}}} 
 
\newcommand{\calo}{{\mathcal O}} 
\newcommand{\oof}[2]{{\mathcal O}_{#1}({#2})} 
\newcommand{\oofp}[2]{{\mathcal O}_{\mathbb{ P}^{#1}}({#2})}

\newcommand{\xel}{(X, L)} 
\newcommand{\sel}{(S, L)}

\newcommand{\Pin}[1]{\mathbb{P}^{#1}}

\newcommand{\lineq}{\sim} 
\newcommand{\gmdp}{ 
\author[G. Besana]{ Gian Mario Besana} 
\address {School of Computer Science, Telecommunications and Information Systems \\DePaul 
University\\ 200 S. Wabash\\Chicago IL 60604 USA} 
 \email{gbesana@cti.depaul.edu}} 
\newcommand{\BSL}{\mathcal{B}(S, L)} 
\newcommand{\BXL}{\mathcal{B}(X, L)} 
\newcommand{\BSV}{\mathcal{B}(S, V)} 
 
\newcommand{\BXV}{\mathcal{B}(X, V)} 
\newcommand{\RSL}{\mathcal{R}(S, L)} 
 
\newcommand{\RXL}{\mathcal{R}(X, L)} 
 
\newcommand{\RXV}{\mathcal{R}(X, V)} 
\newcommand{\Lx}{|L-x|} 
\newcommand{\bslocus}[1]{{\rm Bs} |#1|} 
\newcommand{\bslocusnobars}[1]{{\rm Bs}\, #1 } 
\newcommand{\LA}{\restrict{L}{A}} 
\newcommand{\XV}{(X, V)} 
\newcommand{\phiL}{\varphi_{L}} 
\newcommand{\phiV}{\varphi_{V}} 
\newcommand{\jone}{\mathcal{J}_1(X,L)} 
\newcommand{\jtwo}{\mathcal{J}_2(X,L)} 
\newcommand{\jn}{\mathcal{J}_n(X,L)} 
\newcommand{\ji}{\mathcal{J}_i(X,L)} 
\newcommand{\jonev}{\mathcal{J}_1(X,V)} 
\newcommand{\jtwov}{\mathcal{J}_2(X,V)} 
\newcommand{\jnv}{\mathcal{J}_n(X,V)} 
\newcommand{\jiv}{\mathcal{J}_i(X,V)} 
\newcommand{\jonec}{\mathcal{J}_1(\calL)} 
\newcommand{\jtwoc}{\mathcal{J}_2(\calL)} 
\newcommand{\jthreec}{\mathcal{J}_3(\calL)} 
\newcommand{\jfourc}{\mathcal{J}_4(\calL)} 
\newcommand{\jic}{\mathcal{J}_i(\calL)} 
\newcommand{\jionec}{\mathcal{J}_{i+1}(\calL)} 
\newcommand{\jnc}{\mathcal{J}_n(\calL)} 
\newcommand{\jones}{\mathcal{J}_1(S,L)} 
\newcommand{\jtwos}{\mathcal{J}_2(S,L)} 
 
\newcommand{\calL}{\mathcal{L}} 
\newcommand{\commsquspecial}
{\begin{center}
\begin{picture}(6,9)
\put(1,8){$Y$} \put(5,8){$X$} \put(5,1){$\Gamma$} 
\put(1,1){$\mathbb{F}_2$} \put(1,4){$Z$} 
\put(2,8.2){\vector(1,0){2.5}} \put(2,1.3){\vector(1,0){2.5}} 
\put(1.5,7.5){\vector(0,-1){2}} \put(5.3,7.5){\vector(0,-1){5}} 
\put(1.5,3.5){\vector(0,-1){1.7}} \put(3,8.7){$\sigma$} 
\put(5.7,5){$\pi$} \put(3,1.8){$\nu$} \put(0.5,6.5){$\varphi_2$} 
\put(0.5,2.5){$\varphi_1$} 
\end{picture} 
\end{center}} 
\newcommand{\Pic}{\text{\rm Pic}} 
 
\newtheorem{theo}{THEOREM}[section]
\newtheorem{lemma}[theo]{Lemma}
\newtheorem{cor}[theo]{Corollary}
\newtheorem{prop}[theo]{Proposition}

\theoremstyle{definition}
\newtheorem{dfntn}[theo]{Definition}
\newtheorem{example}[theo]{EXAMPLE}
\newtheorem{conj}[theo]{Conjecture}

\newtheorem{rem}[theo]{Remark}

\theoremstyle{remark}

\begin{document}
\title[Peculiar Loci]{Peculiar loci of ample and spanned line bundles}
\gmdp
\author[S. Di Rocco]{Sandra Di Rocco}
\address{K.T.H.Matematik - S-10044 Stockholm, Sweden.} \email{sandra@math.kth.se}
\author[A.Lanteri]{Antonio Lanteri}
\address{Dipartimento di Matematica ``F. Enriques"
 - Universit\`a degli Studi di Milano - Via Saldini 50 - 20133
Milano, Italy}
\email{lanteri@mat.unimi.it}
 \maketitle 
\setlength{\unitlength}{.5cm} 
%
%
\begin{abstract} The {\it bad locus} and the {\it rude locus} of an ample and base point free linear system 
 on a 
smooth complex projective variety are introduced and studied. 
Polarized surfaces of small degree, or whose degree is the square 
of a prime, with nonempty bad loci are completely classified. 
Several explicit examples are offered to describe the variety of 
behaviors of the two loci. 
\end{abstract} 
 
\section{Introduction} 
%
%
 
Let $X$ be a smooth complex projective variety of dimension $n\geq 2$ and let $L \in \Pic (X)$ be an 
ample line bundle spanned by a vector subspace $V \subseteq H^0(X, L)$. 
 
This work is concerned with the study of two special loci associated with the pair $(X, V)$ (see section 
\ref{notation} for notation). 
\begin{dfntn} 
The {\it bad locus} of $\XV$ is the set 
$$\BXV = \{ x \in X | \text{ every divisor in } |V-x| \text{ is reducible}\}.$$ 
A point $x \in \BXV$ is called a {\it{bad 
point}} for 
$(X,V).$ 
\end{dfntn} 
 
\begin{dfntn} 
\label{rudedef} 
The {\it  rude locus} of $\XV$ is the set 
$$\RXV = \{ x \in X ||V-2x|\neq \emptyset \text{ and every divisor in } |V-2x| \text{ is reducible}\}.$$ 
A point $x \in \RXV$ is called a {\it{rude point}} for $(X,V).$ 
\end{dfntn} 
 
 The first aim of this work is to study $\BXV.$ 
This seems to be a very basic question, but to our knowledge it 
was not previously explicitly discussed in the literature.  The 
phenomenon of bad points turns out, as a consequence of Bertini's 
second theorem, to be intrinsically two-dimensional, see Theorem 
\ref{antotheo}, {\it i).} Moreover the bad locus is a finite set, 
when nonempty, see Corollary \ref{bisfinite}. This last result 
follows from the relationship between the bad locus and the 
$n$-th jumping set of $(X,V),$ as introduced in \cite{la-pa-so2}, 
defined as 
$$\jnv:=\{x \in X \ | \ |V-x| = |V-2x|\}.$$ Indeed 
it is $\BXV \subseteq \jn$, see Theorem \ref{antotheo}, {\it iv)}, 
and $\jn$ is known to be a finite set, \cite{la-pa-so2}, Theorem 
1.2. 
 
 Although in  most cases $\mathcal{B}(X, V)$ is actually empty, for example if $V = H^0(L)$ and 
$L$ is very ample, see Corollary \ref{bisfinite}, or if $\xel$ is 
a scroll, see Corollary \ref{scrollshavenoB}, Section 
\ref{Bexamples} contains concrete examples of surfaces for which 
$\mathcal{B}(X, V) = \jtwov \not=\emptyset.$ These examples show 
that bad points indeed occur and that the result of Theorem 
\ref{antotheo}, {\it iv)} is effective. 
 
Most examples of nonempty bad loci that we present occur for 
surfaces of general type. Nonetheless a non empty $\BXV$ is a very 
rare phenomenon. This suggests the problem of characterizing ample 
and spanned line bundles admitting bad points. An effective, 
complete classification of surfaces with non trivial bad locus and 
degree up to 11 or equal to  $p^2$ where $p$ is prime is 
presented in Section \ref{Bclass}. 
 
 Our interest in {\it bad points} arose from the study of Seshadri constants $\epsilon(\mathcal L,x)$ 
of certain line bundles $\mathcal{ L}$ on a surface $S,$ which are 
slightly more general than ample and spanned ones. In particular 
we were interested in locating the points $x \in S$ where $\epsilon(\mathcal {L},x)$ jumps. For an ample and spanned line 
bundle $L$ on $X$ one has $\epsilon(L,x) \geq 1$ for all $x \in
X$, see for example \cite{Laz2} for $n=2.$ On the other hand one 
can easily see, Proposition \ref{Seshadri2}, that $\epsilon(L,x)
\geq 2$ if $x \in \mathcal B(X,L).$ This leads to the fact that if 
$(X,L)$ is covered by lines then ${\mathcal B}(X,L) = \emptyset,$ 
Corollary \ref{scrollshavenoB}. 
 
The second goal of this work is the study of the rude locus, as in 
Definition \ref{rudedef}. The bad and the rude locus are 
intimately connected. It follows directly from the definitions 
that a bad point is also a rude point, but more precisely 
\begin{equation}
\label{basvjn}
\BXV = \RXV \cap \jnv.
\end{equation}
Notice that the above intersection is nonempty only if $n = 2,$ as 
mentioned above, Theorem \ref{antotheo}, {\it i)}. On the other 
hand there are many rude points that fail to be bad. Indeed, while 
$\BXV$ is a finite set, the behavior of $\RXV$ spans the whole 
gamut of possibilities. In section \ref{Rexamples} examples are 
discussed in which $\RXL$ is respectively empty (Example 
\ref{RDelPezzo2} for $n \ge 3,$ Example \ref{RDelPezzo1} for $n
\ge 4$), a finite set (Example \ref{RDelPezzo2} for $n=2$), a 
divisor with a finite set removed (Example \ref{RS4*}), a divisor 
(Example \ref{RVeronese} for $s=\dim{(|V|)}=2,$ Example 
\ref{RQnO1} for $n=2$), a dense Zariski open subset (Example 
\ref{RDelPezzo1} for $n = 3$), the union of a dense Zariski open 
subset and a finite set (Example \ref{RDelPezzo1} for $n = 2$), 
the whole variety $X.$ 
 The remarkable phenomenon of having $\RXV~=~X$ is addressed in 
Conjecture \ref{congettura}. Notice also that the relationship 
between $\RXV$ and $\BXV$ seen in \brref{basvjn} becomes extreme 
when $\xel$ is a scroll, as it is $\RXL=X$ while, as mentioned 
above, $\BXL = \emptyset.$ 

The authors would like to thank Helmut Epp, Dean of DePaul C.T.I., the Gustafsson Grant 2000-2002 for visiting scientists at K.T.H., the 
NFR 2000-2002 travel grant, and the MURST of the Italian 
Government (in the framework of the PRIN, Cofin 2000, Geometry on 
Algebraic Varieties) for their support. 

The material in this paper is, in part, based upon work supported by the National Science Foundation (NSF) under Grant No. 0125068.  Any opinions, findings and conclusions or recommendations expressed in this material are those of the authors and do not necessarily reflect the views of the NSF.
 
%
%
\section{Notation and Preliminary Results} 
\label{notation} 
 Throughout this article $X$ denotes  a smooth, 
connected, projective variety of dimension $n$ defined over the 
complex field $\mathbb{C}.$ Sometimes we call such an $X$ an 
$n$-fold. If $n=2$ the variety is often denoted with $S.$ Its 
structure sheaf is denoted  by ${\mathcal O}_X$ and the canonical 
sheaf of holomorphic $n$-forms on $X$ is denoted by $K_X$ or 
simply $K$ when the ambient variety is understood. Cartier 
divisors, their associated line bundles and the invertible sheaves 
of their holomorphic sections are used with no distinction. Mostly 
additive notation is used for their group. Given two divisors $L$ 
and $M$ we denote linear equivalence by $L \sim M$ and numerical 
equivalence by $L \equiv M.$ 
 For any coherent sheaf ${\mathcal F}$ on $X$,  $h^i(X,{\mathcal F})$ is the complex dimension of 
$H^i(X,{\mathcal F}).$ When the ambient variety is understood, we often write $H^i({\mathcal F})$ 
and $h^i({\mathcal F})$ respectively for $H^i(X,{\mathcal F})$ and $h^i(X,{\mathcal F}).$ Let $L$ be a  line bundle on $X.$ If $L$ 
is ample, the pair $\xel$ is called a {\it polarized variety}. For 
a subspace $V \subseteq H^0 \xel$ the following notations are 
used:\\ 
$|V|,$ the linear system associated with  $V;$\\ 
$|V-rx|,$ the linear system of divisors in $|V|$ passing trough a point $x \in X$ with multiplicity at least $r;$\\ 
$\bslocus{V},$ the base locus of the linear system $|V|$;\\ 
$\phiV,$ the rational map given by $|V|;$ \\ 
If $V = H^0(L)$ we write $L$ instead of $V$ in all of the above. 
For a line bundle $L \neq \calo_X $ we say that $|L|$ is {\it 
free} if $\bslocus{L}=\emptyset$ or equivalently if $L$ is 
spanned, i.e. generated by its global 
sections;\\ 
The following notation and definitions  are also used:\\ 
$g = g(X, L),$ the {\it sectional genus} of $\xel$, defined by $2g-2=L^{n-1} (K_X+(n-1)L).$ \\ 
$\mathbb{F}_e,$ the Segre-Hirzebruch surface of invariant $e$.\\ 
For the special polarized varieties arising in adjunction theory 
we adopt the usual adjunction theoretic terminology, e.g. see 
\cite{BESO}. In particular, in accordance with it, $(\Pin{1}
\times \Pin{1}, \oof{\Pin{1} \times \Pin{1}}{1,1})$ is not a 
scroll.

 When $S$ is a $\Pin{1}$-bundle over a curve with fundamental section $C_0$ and generic  fiber $f$ we 
have $\text{Num}(S) ={\mathbb Z} \oplus {\mathbb Z},$ generated by the classes of $C_0$ and $f.$ \\ 
Jumping sets of an ample line bundle were introduced in \cite{la-pa-so2}. For the convenience of the 
reader their definition and basic properties are listed below. 
\begin{dfntn}[\cite{la-pa-so2}, Section 1] 
Let $L$ be an ample line bundle on $X,$ spanned by a subspace $V
\subseteq H^0(X, L).$ The set $$\jiv = \{ x \in X | {\rm rk}
(d\phiV)_x \le n -i\}$$ is called the $i-th$ {\it jumping set} of 
$(X,V).$ When $V = H^0(L)$ we set  $\ji = \jiv.$ 
\end{dfntn}
Note that $\jiv$ has a natural scheme-theoretic structure. 
The following chain of inclusions is immediate: 
\begin{equation}
\label{inclusionsJn}
\jonev \supseteq \jtwov \supseteq \dots \supseteq
\jnv.
\end{equation}
Observe that, according to the definition, $\jonev$ is the locus where $\phiV$ ramifies, while,
as mentioned in the introduction,
$$\jnv:=\{x \in X \ | \ |V-x| = |V-2x|\}.$$ 
 
\begin{lemma} 
\label{J2andSing} 
 Let $X$ be a smooth variety of dimension $n$ and 
let $L$ be an ample and spanned line bundle on $X.$ Suppose 
$\phiL: X\to Y$ is an $m:1$ cover of a possibly singular 
$n$-dimensional variety $Y.$ Let $\Delta$ be the branch locus of 
$\phiL.$ Then 
$$ \phiL(\jtwo) \subseteq \left(\text{Sing}(Y) \cup \text{Sing} (\Delta)\right)$$ 
\end{lemma} 
\begin{proof} 
Let $x \in \jtwo \subseteq \jone$ and assume, by contradiction, that  $y =
\phiL(x) \not \in \text{Sing}(Y) \cup \text{Sing} (\Delta).$ The smoothness of
$Y$ and $\Delta$ at $y$ allows us to choose local coordinates 
$(y_1,\ldots, y_n)$ on $Y$ centered at $y$ and $(x_1,\ldots, x_n)$ 
on $X$ centered at $x$ such that in suitable neighborhoods of $x$ 
and $y$ the map $\phiL$ can be written as: 
\begin{gather}
y_1 = x_1^k \notag \\
y_2 = x_2 \notag\\
\ldots \label{localcover} \\
y_n = x_n \notag
\end{gather}
where $2\le k \le m$ and where $y_1 = 0$ is the local equation of 
$\Delta.$ This shows that ${\rm rk} \, (d\phiL)_x = n-1,$ and thus  $x \not \in \jtwo,$ contradiction. \qed \end{proof}

 For the convenience of the reader we also  recall the statement of the 
so-called {\it Bertini's second theorem} and one of its 
corollaries. A nice historical account of the theorem can be found 
in Kleiman's article \cite{kl2}. We sketch the proof of the 
corollary as we do not know of a literature reference for it. 
\begin{theo}[Bertini] 
\label{secondobertini} 
 Let $X$ be a complex projective variety. Let $\Sigma$ be a linear system on $X$ without fixed components, 
  such that every 
$D \in \Sigma$ is reducible. Then $\Sigma$ is composed with a 
pencil $\Lambda.$ 
\end{theo} 
\begin{cor} 
\label{luce} In the hypothesis of  Theorem {\rm 
\ref{secondobertini}}, if the base locus $\bslocusnobars{\Sigma}$ 
is nonempty and finite, then $\dim{X} = 2$ and the pencil 
$\Lambda$ is rational. 
\end{cor} 
\begin{proof} 
Let $p: X --> Y \subseteq \Pin{N} $ be the rational map given by 
$\Sigma.$ The fact that $\Sigma$ is composed with a pencil means 
that $p$ factors through a curve $C,$ say of genus $g,$ i.e., 
$p=\sigma \circ \tau,$ where $\tau : X --> C$ has connected 
fibres, and $\sigma :C \to Y$ is a finite morphism. Any two 
general fibers of $\tau$ meet exclusively along the base locus of 
$\Sigma.$ Because two general elements of a linear system meet in 
codimension one and because the base locus of $\Sigma$ is finite, 
we conclude that $\dim{X} = 2.$ Resolving the indeterminacies of the map $p$ will produce at least one exceptional curve $E$ that dominates $C.$ Thus $g=0.$ 
\qed \end{proof} 
 
\subsection{General Results on Bad Points} 
A reinterpretation of the above second Bertini theorem in our 
context, sheds light on the bad locus of an ample and spanned line 
bundle. 
 
%
%
 
%
%
\begin{theo} 
\label{antotheo} Let $L$ be an ample line bundle on a smooth 
complex projective variety $X$ of dimension $n \geq 2$, spanned by 
a vector subspace $V \subseteq H^0(X, L)$. Let $x \in \BXV.$ Then 
\begin{itemize} 
\item[i)] $n=2;$ \item[ii)]There is an ample line bundle $A$ on 
$X$ with $h^0(A) \ge 2$ such that every $D \in |V-x|$ is of the 
form $D = A_1 + A_2 + \dots A_r,$ $r\ge 2,$ $A_i \in \Lambda,$ 
where $\Lambda \subseteq |A|$ is a linear pencil and thus $L\sim
rA;$ \item[iii)] $|V-x| = |V-rx|$ where $r$ is as in ii); 
\item[iv)] $\mathcal{B}(X,V) \subseteq\mathcal{ J}_n(X,V); $ 
\item[v)] $\phiV^{-1}(\phiV(x)) \subseteq \mathcal {B}(X,V).$ 
\end{itemize} 
\end{theo} 
\begin{proof} 
First of all note that because $L$ is ample and spanned by $V$, 
the base locus of $|V-y|$ is a finite set for all $y\in X,$ and in 
particular $|V-y|$ has no fixed component. As $x$ is a bad point, 
Corollary \ref{luce} can be applied to the linear system $\Sigma =
|V -x|.$ Thus $n=2$ and $|V-x|$ is composed with a rational 
pencil. This means that every $D \in |V-x|$ has the form 
$D=A_1+\dots +A_r$, $r\geq 2$ where $A_i \in \Lambda,$ where 
$\Lambda$ is a linear pencil contained in the complete linear 
system associated with a line bundle $A$ with $h^0(A) \ge 2.$ Thus 
$L \sim rA.$ Moreover $x \in A_i$ for some $i=1, \dots ,r$. 
Because $x$ belongs to all $D \in |V-x|$, necessarily $x$ must 
belong to infinitely many elements $A_i \in \Lambda,$ and 
therefore to all of them and thus $x$ is a point of multiplicity 
$\geq r$ for every $D \in |V-x|$. Thus $|V-x| = |V - rx|.$  In 
particular, since $r \geq 2,$ this shows that $x \in \jnv$ and 
thus $\mathcal{B}(X, V) \subseteq \jnv$. Finally, notice that 
\begin{equation}
\bslocus{V-x}= \{y_1, y_2, \dots , y_s \} = \phiV^{-1}(\phiV(x)),
\label{Antonio1.1.1}
\end{equation}
(where $y_1 = x$). In other words $|V-x|=|V-y_2|= \dots =
|V-y_s|$. In particular this shows that if $x \in
\mathcal{B}(X,V)$, then every $y_i$ is also in $\mathcal{B}(X,
V)$, i.\ e., $$\bslocus{V-x} = \phiV^{-1}(\phiV(x)) \subseteq
\mathcal{B}(X, V).\ \ \ \qed$$  \end{proof}
 
\begin{cor}
\label{bisfinite} Let $L$ be an ample line bundle on a smooth 
complex projective variety $X$ of dimension $n=2,$ spanned by a 
vector subspace $V \subseteq H^0(X, L)$. Then 
\begin{itemize} 
\item[i)] $ \BXV$ is a finite set; 
 \item[ii)] $\BXL = \emptyset$ 
if $L$ is very ample. 
\end{itemize} 
\end{cor} 

\begin{proof} 
Statement {\it i)} follows from the fact that $\dim{(\jiv)} \le
n-i,$ see \cite{la-pa-so2}, Theorem 1.2.  Statement {\it ii)} is 
immediate from the above Theorem \ref{antotheo} and 
\brref{inclusionsJn}. \qed \end{proof} 
 
The distinction between $\BXV$ and $\BXL$ for $V \subset H^0(X,L),$ is subtler than one 
might expect. Clearly $\BXL \subseteq \BXV$ for all nonempty subspaces $V.$ The opposite 
inclusion is in general not true, as the following example shows. 
\begin{example} 
Consider the pair $(S,L)=(\Pin{2}, \oofp{2}{r})$ ($r \geq 2$) and 
let $\Sigma  \subset \Pin{N}$ ($N=\binom{r+2}{2} -1$) be the image 
of $S$ in the embedding $\phiL$. For any point $x \in S$ let $x' =
\phiL(x)$ and let $\text{Osc}^{r-1}_{x'}(\Sigma)$ be the 
$(r-1)$-th osculating projective space to $\Sigma$ at $x'$. Note 
that $\dim{(\text{Osc}^{r-1}_{x'}(\Sigma))}=\binom{r+1}{2} -1$, as 
$L$ is even $r$-jet ample, see \cite{be-sok}. Now let $P$ be a 
general hyperplane of $\text{Osc}^{r-1}_{x'}(\Sigma)$, so that $P$ 
does not intersect $\Sigma$. Note that $N=\dim{(P)} + r+2$. So we 
can consider the projection $\pi_{P}:\Pin{N} - - \to \Pin{r+1}$ 
from $P$ to a $\Pin{r+1}$ skew with $P$. Because $P$ does not 
intersect $\Sigma$, the restricted map $\pi_{P |\Sigma}:\Sigma \to
\Pin{r+1}$ is a morphism. This means that the vector subspace $V
\subset H^0(S,L)$ corresponding to that $\Pin{r+1}$ spans $L$. 
Note that every hyperplane of $\Pin{N}$ containing $P$ and $x'$ 
also contains $\text{Osc}^{r-1}_{x'}(\Sigma)$, hence it is 
osculating of order $r-1$ to $\Sigma$ at $x'$ and therefore the 
section it cuts out on $\Sigma$ has a singular point of 
multiplicity $r$ at $x'$. But such sections correspond to elements 
of the linear system $|V -x|$. So this proves that $|V-x|=|V-rx|$. 
On the other hand the hyperplane sections of $\Sigma$ are 
isomorphic to plane curves of degree $r$. Hence imposing a 
singular point of multiplicity $r$ at $x$ implies reducibility in 
$r$ irreducible components passing through $x$. Therefore $x \in
\BSV$.  Indeed $\BSV= \jtwov =\{x\}$. Note however that $\BSL=
\emptyset$, since $L$ is very ample. Moreover this shows that 
every $x \in S$ can become a bad point for a suitable 
$(r+2)$-dimensional vector subspace $V \subset H^0(S,L)$ spanning 
$L$. 
\end{example} 
 
The following results establish  connections between bad points 
and Seshadri constants. For background material on the latter, see 
\cite{Laz2}, Section 5. 
\begin{prop} 
\label{Seshadri2} 
Let $L$ be an ample and spanned line bundle on an $n$-fold $X$, $n \geq 2$. If $\epsilon(L,x) < 2$, then $x
\not\in \jn$; in particular $x \not\in \BXL$. 
\end{prop} 
\begin{proof} 
 
Since $\epsilon(L,x) < 2$ there is an irreducible curve $C \subset X$ through $x$ such that $LC < 2\
\text{mult}_x(C)$. Because $\bslocus{L-x}$ is a finite set, there is an element 
$D \in |L-x|$ not containing $C$. Thus $$2\ \text{mult}_x(C) > D C \geq \text{mult}_x(D)
 \text{mult}_x(C),$$ which shows that $x$ is a smooth point of $D$. Hence $x \not\in \jn$.
\qed \end{proof}
 
Recall that if $\ell$ is a line of $(X,L)$ then $\epsilon(L,x) = 1$ for every $x \in \ell$ and thus one gets the
following: 
\begin{cor} 
\label{scrollshavenoB} Let $L$ be an ample and spanned line bundle on an $n$-fold $X$, $n \geq 2$. If 
$(X,L)$ is covered by lines then $\jn = \BXL = \emptyset$. In particular this happens for 
scrolls over bases of any dimension $\leq n-1$. 
\end{cor}
\begin{lemma}\label{Seshadrilemma}
Let $L$ be an ample and spanned line bundle on a surface $X,$ 
with $\BXL \ne \emptyset.$ Let $L \sim rA$ be as in Theorem \ref{antotheo} and consider the subset
 ${\mathcal J}_2(A)=\{ x\in X| |A-x|=|A-2x|\}$. 
\begin{itemize} 
\item [i)] If $x \in X \setminus {\mathcal J}_2(A)$  then  
$$r\leq \epsilon(L,x) \le r A^2;$$ 
 \item [ii)]If $A^2=1$ then $X\setminus {\mathcal J}_2(A)=\{x\in X|\epsilon(L,x)=r\}.$ 
\end{itemize} 
\end{lemma} 
 
\begin{proof} 
Recall from Theorem \ref{antotheo} {\it ii)} that $h^0(A) \ge 2.$ Let $x \in X\setminus {\mathcal J}_2(A)$ and let $C$ be any irreducible curve through $x.$ As $x \not \in {\mathcal J}_2(A),$ there exists an element $\hat{A} \in |A - x|$ which is smooth at $x.$ Thus it must be $LC = rAC = r\hat{A}C \ge r \mathrm{mult}_x(C)$ and thus $\epsilon(L,x) \ge r.$ To see the right hand side inequality, choose an element $C\in |A - x|$. Let $C_1$ be the irreducible component of $C$ passing through $x$. Then $\frac{LC_1}{\mathrm{mult}_x(C_1)} \leq rAC = r A^2$ and thus $\epsilon(L,x)\leq rA^2 .$ This proves {\it i)}. 
 
To see {\it ii)} assume $A^2=1$ and $x \in X\setminus {\mathcal J}_2(A)$.  Then $\epsilon(L,x)=r$ from {\it i)}. On the other hand, let $\epsilon(L,x)=r$ and let's choose $C\in|A -x|.$ Then it must be $\mathrm{mult}_x(C)=1$, which implies that $x \in X\setminus {\mathcal J}_2(A).$ 
\end{proof} 

%
%
\begin{lemma} 
\label{lowgenus} Let $\sel$ be a polarized surface with $L$ ample 
and spanned and $\BSL \ne \emptyset.$ Let $A$ be as in Theorem 
{\rm \ref{antotheo}}, {\it ii)}. Then $g(A)\ge 2$ unless $g(A)=
1,$ $L=2A,$ and $(S, A)$ is a Del Pezzo surface with $A^2=1.$ 
\end{lemma} 
\begin{proof} If $g(A)=0,$ then $A$ is very ample (see for example \cite{fu}, (12.1), (5.1)) and therefore $L$ 
is very ample which contradicts $\BSL\ne \emptyset.$ Let $g(A)=1.$ 
As scrolls have empty bad locus, see Corollary 
\ref{scrollshavenoB}, $(S, A)$ must be a Del Pezzo surface, 
according to \cite{fu}, (12.3). If $A^2 \ge 3$ then $A$ is very 
ample, contradiction. Therefore $A^2=1,2.$ If $A^2=2,$ then $2A$ 
is very ample and therefore $L$ is very ample, contradiction. If 
$A^2=1,$ then $3A$ is very ample, therefore it must be $L=2A,$ 
i.e. $(S,L)$ as in the statement. \qed \end{proof}

%
%
\begin{lemma} 
\label{h0L} Let $\sel$ be a polarized surface with $L$ ample and 
spanned and $\BSL\ne \emptyset.$ Let $A$ and $r$ be as in Theorem 
{\rm \ref{antotheo}}, {\it ii)}. Then 
\begin{itemize} 
\item[i)] $h^0(L)\ge r+2,$ in particular $h^0(L)\ge 4;$ 
\item[ii)]Equality holds in  i) if and only if $h^0(A) = 2;$ 
\item[iii)]If equality holds in i) then $\phiL$ is not a 
birational morphism. 
\end{itemize} 
\end{lemma} 
\begin{proof} 
Recall that for two effective line bundles $M$ and $N$ on a smooth 
variety it is $h^0(M+N)\ge h^0(M)+h^0(N)-1,$ see \cite{BESO}, Lemma 1.1.6. Theorem 
\ref{antotheo} gives $L \sim rA,$ and thus one can use the last 
inequality inductively to get $h^0((r-1)A)\ge r,$  as $h^0(A)\ge
2.$ Now consider the sequence 
\begin{equation}
\label{restricttoA} 0\to (r-1)A \to L \to \LA \to 0.
\end{equation}
 Because $L$ is ample and spanned, it must be 
\begin{equation}
\label{dimofim} k=\dim{(Im(H^0(L) \to H^0(\restrict{L}{A})))}\ge 2.
\end{equation}
 Then \brref{restricttoA} combined with \brref{dimofim} 
 gives
$h^0(L)=k+h^0((r-1)A)\ge r+2.$ 
 
Assume equality holds in {\it i)}. Then necessarily equality holds 
in \brref{dimofim} and $h^0(A)=2.$ Now assume $h^0(A)=2$ and let 
$x \in \BSL.$ Then $|L-x|$ is composed with the pencil $|A|,$ 
according to Theorem \ref{antotheo}. Therefore $h^0(L-x) = \dim
{(Sym^r(\mathbb{C}^2))} = r+1.$ Because $L$ is spanned, it is 
$h^0(L) = h^0(L-x) + 1 = r+2.$ 
 
If equality holds in {\it i)} then $h^0(A)=2$ as above. Therefore 
$\restrict{\phiL}{A} : A \to \Pin{1}.$ As $g(A) \ne 0,$ see Lemma 
\ref{lowgenus}, $\restrict{\phiL}{A}$ cannot be birational. 
Because $|A|$ is a pencil then $\phiL$ is not birational. \qed 
\end{proof} 
 
\section{Bad Points Exist} 
\label{Bexamples} As already observed, the bad locus of an ample 
and spanned line bundle is a finite set and it is empty when $n\ge
3.$ It is natural to look for an upper bound for its cardinality. 
Combining  Theorem \ref{antotheo} with the inequality 
$\text{Card}(\jtwo) \leq c_2(J_1L)$ \cite{la-pa-so2}, 
(2.6.1), where $J_1L$ denotes the first jet bundle of $L$, we 
immediately get the upper bound $\text{Card}(\BXL) \leq
c_2(J_1L).$ Note however that this bound is very unsatisfactory as 
easy examples show. Actually for $(X,L) = (\Pin{2}, \oofp{2}{2})$ 
we have $\BXL= \jtwo =\emptyset$, while $c_2(J_1L)= 3$. 
Similarly, let $(X,L)$ be a surface scroll over a smooth curve. Then $\BXL = \jtwo = \emptyset$, according 
to Corollary \ref{scrollshavenoB}, while $c_2(J_1L) = L^2.$ The 
many examples described in this section and the classification 
results contained in Section \ref{Bclass}, suggest a better bound 
of the form $\text{Card}(\BXL) \le \frac{L^2}{4}.$ 
 In most situations $\BXL$ is actually empty, but this is not always the case, as the following example 
inspired by \cite{be-la}, p.\  118-120 shows. 
 
\begin{example} 
\label{coversofcones} 
 Let $\mathbb{F}_e$ be  the Segre-Hirzebruch surface of invariant $e \ge 2.$ For any integer $a \geq 2$ 
consider the line bundle $B = \oof{\mathbb{F}_e}{aC_0+(ae-1)f}$. 
Note that the linear system $|(ae-1)C_0 + e(ae-1)f|$ contains a 
smooth irreducible curve, say $R$ \cite{H}, Corollary 2.18 (b), 
p.\ 380. Let $\varphi:Y \to \mathbb{F}_e$ be a cyclic cover of 
degree $e,$ branched along $C_0 + R \in |eB|$. Note that $Y$ is a 
smooth surface, since $C_0$ and $R$ are smooth and $C_0 R = 0$. 
Let $E: = \varphi^{-1}(C_0)$. Since $C_0$ is a component of the 
branch divisor we have $\varphi^* C_0 = eE$ and then we see from 
$$e^2 E^2 = (\varphi^* C_0)^2 = e C_0^2 = -e^2$$ that $E$ is a 
$(-1)$-curve. Let $\sigma:Y \to X$ be its contraction. Then $X$ is 
a smooth surface; call $x_0 \in X$ the point $\sigma(E)$. Recall 
that $\mathbb{F}_e$ is the desingularization $\nu: \mathbb{F}_e
\to \Gamma_e$ of the cone $\Gamma_e \subset \Pin{e+1}$ over the 
rational normal curve of degree $e$ and that $C_0 + ef =
\nu^*\oof{\Gamma_e}{1}$. Then we get the following commutative 
diagram 
\commsqu{Y}{X}{\Gamma_e}{\mathbb{F}_e}{\sigma}{\pi}{\nu}{\varphi} 
where $\pi$ exhibits $X$ as a cyclic cover of degree $e$ of the 
cone $\Gamma_e$ branched at the vertex $v= \nu(C_0)$ and along the 
transverse intersection with a hypersurface of degree $ae-1$. Note 
also that $\pi^{-1}(v) = \{x_0\}$. 
 
Now, for any integer $b$, with $2 \leq b \leq a$, let $L_b: =\pi^* \oof{\Gamma_e}{b-1}$. Then $L_b$ is an 
ample and spanned line bundle on $X$. For simplicity set $\widetilde L_b = \sigma^* L_b$. Then, due to the 
commutativity of the diagram above, we get 
\begin{equation}
\label{antonio211} \widetilde L_b = \varphi^* (\nu^*
\oof{\Gamma_e}{b-1}) = \varphi^*
\oof{\mathbb{F}_e}{(b-1)(C_0+ef)}.
\end{equation}
 
 Recall that $$\varphi_*(\calo_Y) = \calo_{\mathbb{F}_e} \oplus B^{-1} \oplus \dots \oplus
B^{-(e-1)}$$ (e.\ g.\ see \cite{BPV}, Lemma 17.2, p.\ 43). Then projection formula combined with (2.1.1), 
and the inequality $b \leq a$ give
\begin{xalignat}{1}
 h^0(L_b) &= h^0(\widetilde L_b) \notag\\
                 &= h^0(\varphi_* \widetilde L_b) \notag \\
                 &= h^0((b-1)(C_0+ef)) + h^0((b-1-a)C_0 +(e(b-a)+1-e)f) + \dots \notag \\
                                &= h^0(\nu^*\oof{\Gamma_e}{b-1}) \notag\\
                                &= h^0(\oof{\Gamma_e}{b-1}). \notag
\end{xalignat} 
 
 Thus  \begin{equation}
|L_b| = \pi^* |\oof{\Gamma_e}{b-1}|. \label{antonio212}
\end{equation}

\begin{prop} 
\label{Bforcoversofcones} 
 Let $(X,L_b)$ be as in the current Example. Then $\mathcal{J}_2(X,L_b)=\{x_0\}$ for 
every 
$b$ with $2 \leq b
\leq a$, while
$$ \mathcal{ B}(X, L_b) = \begin{cases} \emptyset &\text{if $b > 2,$}\\ \{x_0\}&\text{if $b=2.$}
\end{cases}
$$
\end{prop}

\begin{proof}The above construction shows that the branch locus of $\pi$ is smooth, 
and therefore Lemma \ref{J2andSing} gives $\mathcal{B}(X, L_b)
\subset \mathcal{J}_2(X, L_b) \subset \{x_0\}.$  By 
\brref{antonio212}, the linear system $|L_b - x_0|$ consists of 
the elements $\pi^*D$ where $D$ is a curve cut out on $\Gamma_e$ 
by a hypersurface of degree $b-1$ passing through $v$. Therefore 
for every $b \geq 2$ all elements in $|L_b - x_0|$ are singular at 
$x_0$, i.\ e., $x_0 = \mathcal {J}_2(X,L_b)$. 
 
Let $b > 2$. Then the general element $D$ as above is irreducible. 
Moreover, if $D$ is in general position with respect to $\nu(R)$, 
then $\pi^* D$ itself is irreducible. This shows that 
$\mathcal{B}(X,L_b) = \emptyset$. On the contrary if $b=2$ then 
all elements $D$ as above consist of $e$ lines  of the ruling of 
$\Gamma_e$ and then $x_0 = \mathcal{B}(X,L_b)$. \qed \end{proof} 
 
Note that $$L_b^2 = e(b-1)^2(C_0+ef)^2 = e^2(b-1)^2,$$ by \brref{antonio211}. In particular $L_b^2=4$ if 
and only if $e=b=2$. Moreover, apart from this case, $L_b^2 \geq 9$, with equality if and only if $e=3, b=2$. 
 
Note also that if $e=a=b=2$ then $X$ is the Del Pezzo surface with $K_X^2=1$ and $L_b=-2K_X$ (e.\ g.\ see 
\cite{BESO}, Example 10.4.3, p.\ 269). In this case the equality $\mathcal{J}_2(X,L_b) = \{x_0\}$ 
was shown in \cite{la-pa-so2}, p.\ 206. 

\begin{rem}\label{seshadriexample}
In the above case,  $b=e=a=2$,  Lemma \ref{Seshadrilemma} gives 
$$\{x\in X\,|\epsilon(L_b,x)=2\}=X\setminus {\mathcal J}_2(-K_X),$$
where $\mathcal{J}_2 (-K_X)$ is the set of double points of the singular elliptic curves, all irreducible, of the pencil $|-K_X|.$
Notice that all elements of the anti-canonical system are smooth at the base point $x_0,$ as $K_X^2=1.$ This implies that $ \{x_0\}={\mathcal B}(X,L_b)\subset X\setminus {\mathcal J}_2(-K_X).$
Moreover 
$${\mathcal J}_2(-K_X)=\{x\in X\,|\epsilon(L_b,x)=1\}.$$ Let $x\in {\mathcal J}_2(-K_X)$ and $C\in |-K_X-x|.$ As $C$ is an irreducible elliptic curve it is $\mathrm{mult}_x(C)=2$ and thus $\epsilon(L_b,x)\leq 1$. On the other hand because $L$ is spanned it is $\epsilon(L_b,x)\geq 1$, see \cite{Laz2}. 

\end{rem}
If $e=b=2, a=3$, then $X$ is a minimal surface of 
general type with $K_X^2=1$, $p_g(X)=2$, studied by Horikawa, \cite{hor2}, and 
$L_b=2K_X$. 
On the other hand $\text{kod}(X)=2$ except for the Del Pezzo surface just mentioned. In fact we have

\begin{prop} 
 Let $X$, $a \geq 2, e \geq 2$ be as in the current example. Then $X$ is a minimal surface of general type 
unless $e=a=2$. 
\end{prop} 
\begin{proof} Recalling the expression of $B$, it is (e.\ g., see \cite{BPV}, Lemma 17.1, p.\ 42) that 
$$K_Y = \varphi^*(K_{\mathbb{F}_e}+(e-1)B) = \varphi^*\oof{\mathbb{F}_e}{cC_0+(ce-1)f}$$ where 
$c = ae -(a+2)$. On the other hand $K_Y = \sigma^*K_X + E$ and $\varphi^*C_0 = eE$. We thus get 
\begin{equation}
e\sigma^*K_X = eK_Y - \varphi^*C_0 = \varphi^*
\oof{\mathbb{F}_e}{\lambda(C_0+ef)}, \label{antonio231}
\end{equation}
 where $\lambda = ce-1$. Note that if $e \geq 3$ then $\lambda >0$, while, if $e=2$ then $\lambda = 2a-5
>0$ unless $a=2$. So, except for the case $e=a=2$, \brref{antonio231} shows that $\sigma^*K_X$ is nef and
big, so being $\oof{\mathbb{F}_e}{C_0+ef}$. This in turn implies 
that $K_X$ itself is nef and big. The nefness shows that $X$ is 
minimal with $\text{kod}(X) \geq 0$ and then the bigness gives 
$\text{kod}(X)=2$. \qed \end{proof}
\end{example}

 By the technique of bidouble covers we can construct another example of some interest in 
itself. 
\begin{example} 
\label{bidoublecovers} 
 Let $\Gamma \subset \Pin{3}$ be the quadric cone, let $\nu: \mathbb{F}_2 \to \Gamma$ be the desingularization, 
 and let $v = \nu(C_0)$ be the vertex of $\Gamma.$ 
 
For any integer $\alpha \geq 1$ consider the line bundle 
$\mathcal{L}=\oof{\mathbb{F}_2}{\alpha(C_0+2f}$. The linear system 
$|2\mathcal{L}|$ contains a smooth divisor $\Delta_1$. Let 
$\psi_1:Z \to \mathbb{F}_2$ be the double cover branched along 
$\Delta_1$. Then $Z$ is a smooth surface; moreover, since 
$\Delta_1 C_0=0$ we have that $\psi_1^*C_0=C_1+C_2$, where 
$C_1,C_2$ are two smooth non-intersecting curves, both isomorphic 
to $C_0$ and exchanged by the involution defined by $\psi_1$. Thus 
the equality 
$$2C_1^2=2C_2^2=C_1^2+C_2^2=(C_1+C_2)^2 = 2C_0^2=-4$$ 
shows that both are $(-2)$-curves on $Z$. For any integer $\beta
\geq 1$ consider on $\mathbb{F}_2$ the linear system 
$|(2\beta-1)(C_0+2f)|$. We can find in it a smooth curve $D$, 
which is transverse to $\Delta_1$ (to see this recall that 
$\oof{\mathbb{F}_2}{(2\beta-1)(C_0+2f)}
=\nu^*\oof{\Gamma}{2\beta-1}$). Thus $R=\psi_1^*D$ is a smooth
curve on $Z$. Note that $DC_0=0$, hence $RC_1=RC_2=0$. Moreover we 
have that $\oof{\mathbb{F}_2}{C_0+D}\in
2\text{Pic}(\mathbb{F}_2)$. So if we consider the line bundle on 
$Z$ given by $B = \psi_1^*\oof{\mathbb{F}_2}{\beta C_0+(2\beta
-1)f},$ we see that $|2B|$ contains the smooth divisor 
$\Delta_2:=\psi_1^*C_0+\psi_1^*D=C_1+C_2+R$. Let $\psi_2:Y \to Z$ 
be the double cover branched along $\Delta_2$. Then $Y$ is a 
smooth surface; moreover since $C_i$ is in the branch divisor for 
$i=1,2$, we have that $\psi_2^*C_i=2E_i$, where $E_i$ is a smooth 
rational curve on $Y$. Thus $4E_i^2=2(C_i)^2=-4$, which means that 
$E_i$ is a $(-1)$-curve on $Y$. Note also that $E_1E_2=0$ since 
$C_1C_2=0$. Let $\sigma:Y \to X$ be the birational morphism 
contracting $E_1$ and $E_2$ and let $x_i=\sigma(E_i)$. Call $\pi:X
\to \Gamma$ the finite morphism of degree $4$ induced by the 
bidouble cover $\psi_1 \circ \psi_2:Y \to \mathbb{F}_2$ via 
$\sigma$ and $\nu$. The following commutative diagram is obtained: 
 
\commsquspecial 
 
 Note that $\pi(x_1)=\pi(x_2)=v$ and $v$ is a branch point of $\pi$, by construction. Now set 
$L_c=\pi^*\oof{\Gamma}{c-1}$, for any integer $c\geq 2$. Then $L_c$ is an ample and spanned line bundle 
on $X$. For shortness set $\widetilde L_c = \sigma^* L_c$. Then, due to the commutativity of the diagram 
above, it follows 
\begin{equation}
\widetilde L_c = \psi_2^*(\psi_1^*
\oof{\mathbb{F}_2}{(c-1)(C_0+2f)}). \label{antonio241}
\end{equation}
Thus, by arguing as in Example \ref{coversofcones}, it is easy to see that if $c
\leq \text{min}\{\alpha,\beta\}$ then 
\begin{equation}
|L_c| = \pi^*|\oof{\Gamma}{c-1}|.
\label{antonio242}
\end{equation}
 
This is the key to prove, in the same way as we did in Proposition \ref{Bforcoversofcones}, that 
$\mathcal{J}_2(X,L_c)=\{x_1, x_2\}$ for every $c$ with $2 \leq c \leq \text{min}\{\alpha, \beta\}$, while $$
\mathcal{B}(X,L_c) =
\begin{cases} \emptyset &\text{if $c > 2,$}\\ \{x_1, x_2\} &\text{if $c=2.$}
\end{cases}$$ 
 
 Finally note that $L_c^2=\widetilde L_c^2 = 4(c-1)^2(C_0+2f)^2=8(c-1)^2$. In particular $L_c^2\geq 8$ 
with equality only if $c=2$. In this case we observe that every element of $|L_c-x_i|=|L_c-x_1-x_2|$ splits in 
two components, which are the images via $\sigma$ of two fibres of $Y \to \Pin{^1}$ belonging to the 
pencil $|\psi_2^*(\psi_1^* f)|$. 
\end{example} 
 
\section{Polarized Surfaces of low degree with non trivial bad loci} 
\label{Bclass} 
%
%
 
 The examples studied in the previous section and Lemma 
\ref{J2andSing} suggest that non empty bad loci commonly occur on 
surfaces that are covers of cones. The fact that $L$ is divisible 
in $\text{Pic}(S),$ as shown in Theorem \ref{antotheo}, imposes 
simple stringent conditions on the possible degree of a polarized 
surface admitting a non trivial bad locus. In particular it 
immediately follows from Theorem \ref{antotheo} that $L^2 = 4,8,9$ 
or $L^2 \ge 12.$  In what follows, polarized surfaces with non 
trivial bad loci of degree up to 9 or equal to the square of a 
prime $p$ are completely classified. Indeed they turn out 
to be finite covers of cones, as anticipated in the examples of 
the previous section. 
\begin{prop} 
\label{B49} Let $(S, L)$ be  a polarized surface with $L$ ample 
and spanned. Assume $\BSL \ne \emptyset.$ Let $L^2 = p^2$ where 
$p$ is prime. Then \begin{itemize} 
 \item[i)] $\BSL = \{x\};$ 
\item[ii)]$\phiL: S \to Y_p$ expresses $S$ as a $p$-uple cover of 
a cone $Y_p\subset \Pin{p+1}$ over a rational normal curve of 
degree $p,$ whose branch locus contains the vertex $v;$ 
\item[iii)]$\phiL(x) = v;$ 
 \item[iv)]In particular, if $p=2$ then $L=2A$ where $(S,A)$ is a Del Pezzo 
 surface of degree one. 
 \end{itemize} 
\end{prop} 
\begin{proof} 
 Let $x \in \BSL.$ Theorem \ref{antotheo} gives $L=pA$ with $A^2=1.$ This implies that every pair of divisors in the 
 pencil $\Lambda \subseteq |A|$ meet 
only at $x.$  As $ |L-x|$  does not have a fixed component, it is 
set theoretically 
$\bslocusnobars{\Lx}=\bslocusnobars{\Lambda}=\{x\}.$ Bertini's 
theorem then gives that  the generic $A_i\in \Lambda$ is 
everywhere smooth. Let us consider a generic smooth $A \in
\Lambda.$  If $g(A) =1$ Lemma \ref{lowgenus} implies that $p=2$ 
and $A = -K_S.$  It is well known and easily seen that in this 
case $\varphi_{-2K_S}$ expresses $S$ as a double cover of a 
quadric cone in $\Pin{3},$ see for example \cite{BESO},
Ex.\,10.4.3, p.\ 269. Therefore $g(A)\ge 2,$ and as $A^2=1$ the 
sequence $0 \to \calo_S \to A \to \restrict{A}{A} \to 0$ gives 
$h^0(A)= 2.$ Lemma \ref{h0L} then gives $h^0(L)= p + 2,$ and 
$\phiL$ not birational. In particular $\deg{(\phiL(S))}\ge 2.$ 
From\begin{equation} \label{degdeg} p^2=\deg{(\phiL)}
\deg{(\phiL(S))}
\end{equation} we see that $\phiL$ must be a $p$-uple cover of a surface $Y_p\subset \Pin{p+1}$ of degree $p.$ 
 As $Y_p$ is covered by lines $\phiL(A_i)$ for $A_i \in
|A|,$ all of which meet at one point, $\phiL(x),$ $Y_p$ must be a 
cone with vertex $v =\phiL(x),$ over a rational normal curve of 
degree $p.$  Because $v$ is the only point on $Y_p$ through which 
there are infinitely many lines, $\BSL={x}.$ When $p=2$ it follows 
from \cite{fuhy}, Section 4 that $\sel$ is as in {\it iv)}. \qed 
\end{proof} 
 
%
%
\begin{prop} Let $\sel$ be a smooth surface polarized with an ample and spanned line bundle with $L^2=8.$ 
Assume that $\mathcal{B}(S,L)\ne\emptyset.$ Then 
\begin{itemize} 
\item[i)]$\mathcal{B}(S,L)=\{x_1,x_2\}$ with $x_1\ne x_2;$ 
\item[ii)] $\phiL: S \to Y_2$ expresses $S$ as a quadruple cover 
of a quadric cone $Y_2 \subset \Pin{3}$ with vertex $v;$ 
\item[iii)] $\phiL(x_1)=\phiL(x_2)= v.$ 
\end{itemize} 
\end{prop} 
\begin{proof} Let $x\in \mathcal{B}(S, L)$ and let $C$ be a generic element in $\Lx.$ 
Note that $\text{mult}_x(C)\ge 3$ would imply $L^2\ge 9,$ and thus 
it is $\text{mult}_x(C)=2.$  If there exists $y \in
\bslocusnobars{\Lx},$ $y\ne x,$ Theorem \ref{antotheo} implies
that $y$ is also a bad point for $L$ and thus $\text{mult}_y(C)\ge 2.$ These facts, combined with the second Bertini Theorem, imply 
that the only two possible configurations for 
$\bslocusnobars{\Lx}$ are as follows: 
\begin{itemize} 
\item[a)]$\bslocusnobars{\Lx}=\{x\} $ with the intersection index at $x$ $(C_1,C_2)_x=L^2=8$ for two general 
elements of $\Lx.$ 
\item[b)]$\bslocusnobars{\Lx}=\{x_1,x_2\}$ with $x_1\ne x_2,$ and $(C_1,C_2)_{x_1}=(C_1,C_2)_{x_2}=4$ for all 
$C_i \in \Lx.$ 
\end{itemize} {\it \bf Claim} Case a) does not happen. 
 
In case a), the general $C\in \Lx$ would be smooth away from $x$ 
and two general such curves would intersect only at $x.$  This 
fact, combined with the two Bertini theorems, forces the general 
$C$ to be reducible at $x$ with $C=A+B,$ $A$ and $B$ both smooth 
at $x$ and linearly equivalent, $A^2=AB=B^2=2.$   This  implies 
that the tangent lines to $A$ and $B$ at $x$ must be the same and 
that this line does  not change when we vary $C$ generically in 
$\Lx.$ Let now $C_i \in \Lx,$ $i=1,2$  be two general curves, and 
consider the pencil $\mathcal{P}=<C_1,C_2>.$ Choosing suitable 
local coordinates $z_1,z_2$ centered at $x,$ such that the common 
tangent at $x$ has equation $z_1=0,$ the local equations of the 
$C_i$ are $z_1^2 + F_i(z_1,z_2)= 0$ where the $F_i$'s  are 
polynomials of degree at least $3.$ Therefore $\mathcal{P}$ 
contains the curve $C_1-C_2$ which has at least a triple point at 
$x.$ As $L^2=8,$ there can be only one such curve in $\Lx.$ This 
means that $\Lx=\mathcal{P}.$ Indeed, if $\mathcal{P} \subset
\Lx$, we could find a general $C_3 \in \Lx, C_3 \not \in
\mathcal{P}.$ Now with $C_1$ and $C_3$ 
 we could repeat the above  pencil argument and produce a second curve with a triple point, different 
from the first one, which is a contradiction. Therefore $h^0(L-x)
= 2$ and, because $L$ is ample and spanned, $h^0(L)=3,$ but this 
contradicts Lemma \ref{h0L}, {\it i)}. This concludes the proof of 
the Claim. 
 
Let now $\Lx$ be as in b). Again Bertini's theorems give $L=2A $ 
where $A \in  \Lambda,$ a pencil in $|A|,$ and that the generic $A
\in \Lambda$ is smooth  everywhere and goes through both $x_i$'s. 
Notice that $\{x_1,x_2\}\subseteq \BSL\subseteq \jtwos \subseteq
\jones$ and thus the $x_i$'s are both ramification points for the 
finite map $\phiL.$ Because $\phiL(x_1)=\phiL(x_2),$ it follows 
that $\deg{( \phiL)}\ge 4.$ As $L^2=8,$  $\deg{(\phiL(S))}\le 2$, 
and thus $h^0(L)\le 4.$ Therefore Lemma \ref{h0L} gives $h^0(L)=4,
h^0(A)=2$ and $\phiL$ not birational. Thus the sequence $0\to A\to
L\to\LA\to 0$ gives $\dim{(Im(H^0(L) \to H^0(\LA)))}= 2.$ 
Therefore $\restrict{\phiL}{A}(A)$ expresses $A$ as a $4:1$ cover 
of  a line. As the pencil $|A|$ sweeps out  the whole $S$, it must 
be $\deg{(\phiL)}\ge 4$ and thus $\deg{(\phiL(S))}= 2.$ Let 
$\mathbb{Q}=\phiL(S).$ As $\mathbb{Q}$ is swept by a pencil of 
concurrent lines, images of the pencil of curves $|A|,$ 
$\mathbb{Q}$ must be  a quadric cone with vertex $v=\phiL(x_i).$ 
Because $v$ is the only point on $\mathbb{Q}$ through which there 
are infinitely many lines, $\BSL=\{x_1,x_2\}.$ \qed \end{proof} 
 
\section{The rude locus} 
\label{Rexamples} The following set of examples is aimed at 
showing how the behavior of the rude locus of an ample and free 
linear system covers a very wide spectrum of possibilities. As 
mentioned in the introduction, $\RXV$ can be empty, a finite set, 
a divisor minus a finite set, a divisor, a dense Zariski open subset, 
the union of a dense Zariski open subset and a finite set or the whole 
variety. 
 
 Let $\xel$ be a polarized 
$n$-fold with $V\subseteq H^0(L)$ a subspace that spans $L.$ In 
all the following examples we will denote with $s$ the projective 
dimension of the linear system $|V|.$ Notice that our assumptions 
on $L$ and $V$ imply that $s=\dim{(|V|)} \ge n.$ As customary we 
set the dimension of the empty set equal to $-1.$ 
\begin{lemma} 
\label{introtoR} Let $\xel$ be a polarized $n-fold,$ and 
$V\subseteq H^0(L)$ be a subspace that spans $L.$ Let $s =
\dim{(|V|)},$ let $x \in X$ and let $\phiV$ be the map associated 
with $|V|.$ Then 
\begin{itemize}
\item[i)]$\dim{( |V - 2x|)} = s - 1 - {\rm rk} (d\phiV)_x$;
 \item[ii)] $\dim{(|V - 2x|)} \ge s-(n+1)$ with equality
holding if and only if  $x \in X \setminus \jonev$;
\item[iii)] If $\dim{(|V|)} = n$ then $\RXV \subseteq
\jonev$ and $|V - 2x| \neq \emptyset$ for all $x \in
\jonev.$
\end{itemize}
\end{lemma}
\begin{proof}Let $\mathfrak{m}_x$ be the maximal ideal of $\calo_{X,x}$ and consider the homomorphism
 $$j_{1,x}: V \to L \tensor \calo_X/\mathfrak{m}_x^2 $$
 sending every section $s \in V$ to its first jet at $x,$ i.e. $j_{1,x}(s) = (s(x), ds(x))$ in a local chart around $x.$
  Because $V$ spans $L$ at $x$ it is ${\rm rk} (j_{1,x}) = 1 + {\rm rk}
  (d\phiV)_x.$
Then {\it i)} follows noting that $|V - 2x| =\mathbb{P}(Ker(j_{1,x})).$ To see {\it ii)} notice that ${\rm rk}
\, (d\phiV)_x \le n$ with equality if and only if $ x \in X
\setminus \jone.$ Assuming $s = n,$ {\it ii)} gives {\it iii)}.
\qed \end{proof}

%
%
\begin{example}[The Veronese surface] 
\label{RVeronese} Let $\xel = (\Pin{2}, \oofp{2}{2})$ and let $V
\subseteq H^0(\Pin{2}, \oofp{2}{2})$ be  a  subspace that spans 
$L.$ Let $s = \dim{(|V|)}.$ Then it is 
$$\RXV = \begin{cases} 
X \text{ \ \ \ \ \ \ \ \ \ \ \  if\  } s \ge 3,\\ 
\jonev \text{  \ \  if \ } s = 2. 
\end{cases}$$ 
Assume first that $s \ge 3.$ Because every singular conic is 
necessarily reducible, it is enough to show that for every $x \in
X,$ $|V- 2x|$ is not empty. This follows from the condition $s \ge
3$ and Lemma \ref{introtoR},  {\it ii)}. 
 
Now assume $s = 2.$ Because of Lemma $\ref{introtoR}$ we need to 
show that if $x \in \jonev$ then $x \in \RXV.$  It is enough to 
show that $|V - 2x|$ is nonempty. This follows again from Lemma 
\ref{introtoR}, {\it ii)}. 
\end{example} 
 
 The fact that being singular and being reducible are equivalent for rational curves  is responsible for the 
peculiarity of the two dimensional case among the polarized 
varieties $(\Pin{n}, \oofp{n}{2}).$ The following example 
describes the general picture. 
 
%
%
\begin{example} 
\label{RPnO2} Let $(X,L)=(\Pin{n}, \oofp{n}{2})$, with $n \geq 2$ 
and let $V \subseteq H^0(X,L)$ be a general subspace, with $s 
=\dim{(|V|)} \geq n+1.$ Let $\mathcal{W} \subset |L|$ be the 
subvariety parameterizing all reducible quadric hypersurfaces of 
$\Pin{n}$. 
 
Notice that  $\mathcal{W}$ is isomorphic to the second symmetric 
power of the dual of $\Pin{n}$, hence $\dim{(\mathcal{W})} = 2n$. If 
$\mathcal{W} \cap |V|=\emptyset$, then 
$\mathcal{R}(X,V)=\emptyset$. Now suppose that $\mathcal{W}$ meets 
$|V|$, and inside $X \times (|V|\cap \mathcal{W})$ consider the 
following incidence variety: $$\mathcal{I}:=\{(x,D)| x \in 
\text{Sing}(D)\}.$$ Let $p$ and $q$ be the morphisms induced 
on $\mathcal{ I}$ by the projections of $X \times (|V|\cap
\mathcal{W})$ onto the factors. Notice that 
$$p^{-1}(x) = \{D \in |V| \cap \mathcal{W}\ | \ x \in 
\text{Sing}(D)\} = |V-2x| \cap \mathcal{W},$$ and therefore $x \in
\mathcal{R} (X,V)$ if and only if 
\begin{equation}
\label{lafibradip} p^{-1}(x) = |V - 2x|. 
\end{equation} 
 
We claim that \brref{lafibradip} cannot occur for $n \ge 3,$ and 
for a general $x \in X.$ 
 
Since $V$ is general, the codimension of $(|V|\cap \mathcal{W})$ 
in $|V|$ is the same as that of $\mathcal{W}$ in $|L|$, hence 
$$\dim{(|V|\cap \mathcal{W})} = s-\left(\binom{n+2}{2} -1 - 2n\right)=s- \binom{n}{2}.$$ 
 As the singular locus of 
an element $D \in |V|$ broken into two distinct hyperplanes is a 
linear space of codimension $2$, the general fiber of $q$ has 
dimension $n-2.$ Therefore $\dim{(\mathcal{I})} =
s-\binom{n}{2}+n-2$ and then for the general $x \in X$ we have 
that $\dim{( p^{-1}(x))} = s - \frac{1}{2}(n^2-n+4)$. \

Therefore Lemma \ref{introtoR} gives for a general $x$: 
\begin{gather} \dim{(|V- 2x|)} - \dim{(p^{-1}(x))} \ge s - n -1 -(s
- \frac{1}{2}(n^2-n+4)) \\ \notag  = \frac{1}{2} (n - 1) (n - 2)
\end{gather} 
which proves the claim.

This discussion proves that $\mathcal{R}(X,V)$ is contained in 
(or is equal to) a closed Zariski subset of $X.$ It is easy to 
see that $\RXL = \emptyset.$ 
\end{example} 
 
The following example shows that the inclusion given in Lemma 
\ref{introtoR} { \it iii)} can be proper. Notice again the 
striking difference between the behavior of $\RXV$ in dimension 
two versus the higher dimensional cases. 
 
%
%
\begin{example}[Del Pezzo manifolds of degree $2$] 
\label{RDelPezzo2} Let $\xel$ be a Del Pezzo manifold of dimension 
$n \ge 2$, i.e. $-K_X = (n-1)L$ with $L$ ample, $L^n = 2.$ Then 
$L$ is spanned, $h^0(L) = n+1$ and $\phiL : X~\to~\Pin{n}$ is a 
double cover ramified over a smooth quartic hypersurface $\Delta.$ 
For $n = 2,$ recall that a smooth plane quartic admits $28$ 
double-tangent lines $\ell_k.$ Let $\Delta \cap \ell_k = \{
y^k_1,y^k_2\}$ and consider the $56$ points $x^k_j =
\phiL^{-1}(y^k_j)$, $j=1,2$, $k = 1, \dots, 28.$ If $V$ is a 
subspace that spans $L,$ necessarily $|V| =|L|$ because 
$\dim{(|L|)} = n.$ 
 Then $$ 
\RXL = \begin{cases} 
\emptyset \text {\ \ \ \ \ \ if } n\ge 3, \\ 
 \{x_j^k|k=1, \dots,28,\ \ 
j=1,2\} \text{ \ \ \ \ \  if  } n=2. 
\end{cases} 
$$ 
 
Assume first that $n \ge 3.$ Then $\text{Pic}(X) \simeq 
\text{Pic}(\Pin{n}) = \mathbb{Z}$ generated by $L,$ see for 
example \cite{Laz3}, Prop. 3.1. If $x \in \RXL$ and $D \in |L-2x|$ 
then it must be $D = A + B$ for two effective divisors $A$ and 
$B.$ Thus $A=aL$ and $B=bL$ where $a,b \ge 1.$  This gives the 
following chain of equalities in $\Pic (X)$ $$L=D=A+B=(a+b)L$$ 
which is a contradiction as $a+b \ge 2.$ 
 
Let now $n=2.$ Lemma \ref{introtoR} gives $\RXL \subseteq \jone.$ 
For every point $x \in \jone$ the elements in $|L-2x|$ are the 
preimages via $\phiL$ of the lines in $\Pin{2}$ which are tangent 
to $\Delta$ at $y=\phiL(x).$ For a generic $x,$ such a tangent 
line is tangent to $\Delta$ only at $y,$ and therefore the 
corresponding element in $|L-2x|$ is an irreducible curve with 
arithmetic genus one, with one double point. However $|L-2x^k_j|$ 
consists of a curve, $\phiL^{-1}(\ell_k),$  of arithmetic genus 
one  with two double points $x^k_1, x^k_2$ and therefore 
reducible. Thus $\RXL = \{x_j^k|k=1, \dots,28,\ \ j=1,2\}.$ Notice 
that here the inclusion $\RXL \subset \jone$ is strict. 
\end{example} 
 
%
%
 
\begin{example}[Del Pezzo manifolds of degree $1$] 
\label{RDelPezzo1} Let $(X,\mathcal{L})$ be a Del Pezzo manifold 
of degree one, i.e. 
 $-K_{X}=(n-1)\mathcal{L}$, with $\mathcal{L}$ ample and $\mathcal{L}^{n}=1.$ Then $h^0(X,\mathcal{L})=n.$ 
Moreover $\mathcal{L}$ is not spanned and $\text{Bs}|\mathcal{L}|$ 
consists of a single point $x^*$.  Set $L = 2\mathcal{L}$. The 
line bundle $L$ is ample and spanned, with 
$h^{0}(X,L)={{n+1}\choose 2}+1$  and it defines a $(2:1)$ cover, 
$\phiL:X\to \Gamma$ of a cone over the $2$-Veronese manifold of 
 dimension $n-1$, $(\Pin{n-1}, \oofp{n-1}{2})$, where $\phiL(x^{*})= v,$ the vertex of the cone. 
For general results on Del Pezzo manifolds see \cite{fu}, Chapter 
I.

  Even  though the line bundle $\mathcal{L}$ is not spanned, we 
  set 
$$\jic = \{x \in X \setminus\{ 
x^{*}\}| {\rm rk}\,(d\varphi_{\mathcal{L}})_{x}\leq n-i\}.$$ 
\begin{rem} 
\label{jstortopersup} If $n = 2$, as ${\rm rk}(\phiL)_x \le 1,$ it 
is $\jonec = X \setminus{x^*}.$ Moreover $\jtwoc$ is given by the 
finite set of the singularities $x_1, \dots, x_\delta$ of the 
singular curves in the pencil $|\mathcal{L}|.$ It is well known 
that if $X$ is general in the moduli space of Del Pezzo surfaces 
of degree $1$ then $|\mathcal{L}|$ contains exactly $12$ singular 
curves, each of which has a single ordinary node, hence $\delta = 
12.$ 
\end{rem} 
 
\begin{lemma} 
\label{dpzlemma1} 
 
$$\jionec \subseteq \ji \subseteq \jic \cup \{x^*\}.$$ 
\end{lemma} 
\begin{proof} 
First notice that $x^* \in \jn \subseteq \ji$ for all $1\le i\le 
n$  because every $D \in |L - x^*|$ is of the form $D = 
\phiL^*(H)$ where $H$ is a hyperplane section of $\Gamma$ through 
$v$ and thus $D \in |L - 2x^*|.$ 
 
Let now $x \neq x^*$. Let us fix a basis for $H^0(X,\calL)$, 
$\mathfrak{B}=\{s_0,\dots,s_{n-1}\}$. Let $\mathfrak{B}_1 $\ be 
the basis of $Sym^{2}(H^0(X,\mathcal{L}))$ constructed from 
$\mathfrak{B}.$ We can choose an appropriate section $h$ such that 
$\mathfrak{B}_1 \cup \{h\}$ is a basis for $H^0(X, L).$ 
 
Assume ${\rm rk}\,(d\varphi_{\calL})_x = t.$ Then, after choosing 
local coordinates around $x,$ say $x_1, \dots x_n$, and 
renumbering the sections in $\mathfrak{B}$ we can assume that in a 
neighborhood of $x$ it is $s_0 =1,$ $s_i = x_i + \text{\it{higher 
order terms}}$ for $i = 1, \dots, t$ and $s_j \in H^0(\calL - 2x)$ 
for $j = t+1, \dots, n-1.$ From our construction of the basis for 
$H^0(X,L)$ it follows that the only linearly independent rows in 
the matrix $(d\phiL)_x$ are given by $(d(s_0s_i))_x,$ for  $i=1 
\dots t,$ and possibly $(dh)_x.$ Thus ${\rm rk} \, (d\phiL)_x = t$ 
or $t+1.$ The above argument shows that 
 $$  {\rm rk} \, (d\phiL)_x \le{\rm rk} \, (d \varphi_{\calL})_x +1$$ 
from which we deduce that $\jionec \subseteq \ji$. 
 
Assume now that ${\rm rk}\,(d\varphi_{L})_x = t$. Then after 
choosing appropriate local coordinates, $y_1,...,y_n$, around $x$, 
we can assume that there are $t+1$ global sections $h_0,...,h_t\in
\mathfrak{B}_1 \cup \{h\}$  that locally can be expressed as $h_0 
=1,$  $h_k = y_k + \text{\it{higher order terms}}$ for $k = 1, 
\dots, t$. If $h_k\neq h$ for all $k=0,...,t$ then 
$h_k\in\mathfrak{B}_1$ and thus $h_k=s_{i_k}s_{j_k}$, with 
$s_{i_k},s_{j_k}\in\mathfrak{B} $. The expansions of the $h_k's$ 
in local coordinates show that there are $t+1$ sections in 
$\mathfrak{B}$, $s_{i_0},...,s_{i_t}$ such that 
$h_0=s_{i_0}s_{i_0}$ and $h_k=s_{i_0}s_{i_k}$. Locally it is 
$s_{i_0} =1,$  $s_{i_k} = y_k + \text{\it{higher order terms}}$ 
for $k = 1, \dots, t$. Those sections form a non zero $t\times t$ 
minor for the matrix $(d\varphi_{\calL})_x$. It follows that ${\rm 
rk} \, (d\varphi_{\calL})_x = t$. If there is a $k$ such that 
$h_k=h$, then, reasoning as above,  ${\rm rk} \, 
(d\varphi_{\calL})_x = t$ or ${\rm rk} \, (d\varphi_{\calL})_x = 
t-1$. The above argument shows that 
 $$  {\rm rk} \, (d\phiL)_x -1 \le{\rm rk} \, (d \varphi_{\calL})_x\le{\rm rk} \, (d\phiL)_x  $$ 
from which we deduce that $\ji \subseteq \jic \cup \{x^*\}$. \qed 
\end{proof} 
 
\begin{lemma} 
\label{dpzlemman=2} Let $\xel$ be as above with $n=2.$ Then 
$$ \RXL = (X \setminus \jone) \cup \jtwoc \cup \{ x^*\}$$ 
where $\jtwoc$ is as in Remark \rm{\ref{jstortopersup}}. 
\end{lemma} 
\begin{proof}  In dimension $2$, 
 $\phiL: X \to \Gamma$ is a double cover of a quadric 
cone in $\Pin{3}$, ramified along the vertex and $\Delta$, the 
intersection of a smooth hypersurface of degree $3$ with $\Gamma.$ 
The jumping loci are worked out in \cite{la-pa-so2} where it is 
shown that $\jtwo =\{x^{*}\}$ and $\jone = \mathfrak{R} \cup \{ 
x^*\}$, where $\mathfrak{R}$ is the ramification divisor 
$\phiL^{-1}(\Delta).$ If $x\not\in \jone$ then $\dim{(|L-2x|)}=0$ 
i.e. there is a unique singular divisor in $|L - x|$, which must 
be of the form $2D$, where $D$ is the unique divisor in 
$|\mathcal{L}-x|$. Hence $X \setminus \jone \subset \RXL$. We have 
already seen that $x^*\in \BXL \subset\RXL.$ When $x^*\neq x\in 
\jone$ then $\dim{( |L-2x|)}=1$, i.e. there is a pencil of 
singular sections in $|L - x|.$  A section  $D\in |L-2x|$ is given 
by $\phiL^*H$, where $H$ is a hyperplane section of the cone 
$\Gamma,$ tangent to $\phiL(\mathfrak{R}) = \Delta$ at $\phiL(x)$, 
i.e. given by an element of the pencil of planes containing the 
tangent line to $\Delta$ at $\phiL(x)$. The generic such $H$ is 
irreducible, unless the tangent line $\ell_x$ to $\Delta$ at 
$\phiL(x)$ is a line of the cone through $v$, which corresponds to 
the unique element  $T \in |\mathcal{L} -x|.$ This happens only if 
$x \in \jtwoc.$ Notice that in this case every element of the 
pencil of singular sections is reducible as $T+\hat{T}$, with 
$\hat{T}$ varying in $|\mathcal{L} |$ and thus $\jtwoc \subseteq 
\RXL.$ \qed \end{proof}

\begin{lemma} 
\label{dpzlemma2} Let $n\ge 3$ and let $x  \in \RXL,$ then 
$$\dim{(|L-2x|)} = \begin{cases} 
2n -4 \text{ \ \ \ if   } x \not \in \jnc \cup \{ x^*\},\\ 
2n -3 \text{ \ \ \ if   } x \in \jnc,\\ 
2n - 2 \text{ \ \ \ if   } x = x^*. 
\end{cases}$$ 
\end{lemma} 
\begin{proof} 
Consider the following two families of divisors on $X$ 
$$\mathcal{F} = \{ D_1 + D_2 |\  D_i \in |\calL -x|\  i=1,2\}$$ and 
$$\mathcal{G} = \{ D_1 + D_2 |\  D_1 \in |\calL -2x|, D_2 \in |\calL| 
\}.$$ There is a natural $2:1$ cover $|\calL - x| \times |\calL - 
x| \to \mathcal{F}$ and thus 
\begin{equation} 
\label{dimF} 
\dim{(\mathcal{F})} 
=\begin{cases} 2n - 4 \text{ \ \ \ if } x \neq x^*,\\ 
2n -2 \text{ \ \ \ if } x = x^*. 
\end{cases} 
\end{equation} 
Furthermore, $|\calL-2x| \times (|\calL|\setminus|\calL-x|)= 
|\calL-2x| \times \mathbb{C}^{n-1}$ is a dense Zariski open subset 
of $\mathcal{G}.$ Thus 
\begin{equation} \label{dimG} \dim{(\mathcal{G})} = (n-1) + 
\dim{(|\calL-2x|)}. 
\end{equation} 
 Moreover it is 
\begin{equation} 
\label{dimcalL2x} \dim{(|\calL - 2x|)} = \begin{cases} n-2 \text{ 
\ \ \ if } x \in \jnc,\\ 
\le n-3 \text{ \ \ \  if } x \in (X \setminus \jnc) \cup \{ x^* 
\}. 
\end{cases} 
\end{equation} 
Because $\text{Pic}(X) = \mathbb{Z}$ generated by $\calL$, a 
divisor $D \in |L -2x|$ is reducible if and only if $D \in 
\mathcal{F} \cup \mathcal{G}.$ If $x \in \RXL$ then $\dim{(|L - 
2x|)} = \text{max} \{\dim{(\mathcal{F})}, \dim{(\mathcal{G})}\}.$ 
Putting together \brref{dimF}, \brref{dimG}, and \brref{dimcalL2x} 
one obtains the statement. \qed \end{proof}

\begin{prop} 
Let $(X,\mathcal{L})$ be a Del Pezzo manifold  of dimension $n$ 
with $\calL^n = 1.$ Let $L = 2 \calL.$ Then 
$$\RXL = \begin{cases} 
\emptyset \text{ \ \ \ for } n \ge 4,\\ 
X \setminus \jone \text { \ \ \ for } n = 3,\\ 
(X \setminus \jone) \cup \jtwoc \cup \{ x^*\}\text { \ \ \ for } n 
= 2. 
\end{cases} 
$$ 
\end{prop} 
\begin{proof} 
Lemma \ref{dpzlemman=2} gives the statement for $n=2.$ Combining 
Lemma \ref{dpzlemma2}, Lemma \ref{dpzlemman=2} and Lemma 
\ref{introtoR} { \it i)}, it follows that $\RXL = \emptyset$ if $n 
\ge 5$ 
 and $ x^* \in \RXL$ if and only if $ n = 2.$ 
 If $n=4$ Lemma \ref{dpzlemma2} implies that if $x \in \RXL$ then 
 $ x \in \jfourc \cap (X \setminus \jone)$ but this intersection 
 is empty according to Lemma \ref{dpzlemma1}. 
 If $n=3$ and $x \in \RXL$ then Lemma \ref{dpzlemma2} gives $$x \in 
 [\jthreec \cap (\jone \setminus \jtwo)] \cup [(X\setminus 
 \jthreec)\cap (X\setminus \jone)].$$ Lemma \ref{dpzlemma1} gives $$\jthreec \cap (\jone \setminus 
 \jtwo)= \emptyset$$ and $$(X\setminus 
 \jthreec)\cap (X\setminus \jone) = X\setminus \jone.$$ 
 Thus $\RXL \subseteq X \setminus \jone.$ 
 Following the proof of Lemma \ref{dpzlemma2}, and using the same notation introduced there, notice that 
 when $n=3$ it is $|\calL - x| \simeq \Pin{1}$ and thus $\mathcal{F} \simeq 
 \Pin{2}.$ If $x \in X \setminus \jone$ then $\dim{(|L-2x|)} = 2$ 
 and thus $|L- 2x| = \mathcal{F}$ i.e. $x \in \RXL.$ 
\qed \end{proof} 
\end{example} 
 
%
%
\begin{example}[$\mathbb{Q}^n, \oof{\mathbb{Q}^n}{1}$] 
\label{RQnO1} Let $X = \mathbb{Q}^n$ and $L= 
\oof{\mathbb{Q}^n}{1}.$ Let $V \subset H^0(X, L)$ be a spanning 
subspace. Then either $|V| = |L|$ or ${\rm codim}_{|L|}(|V|) = 1.$ 
In the first case if $n=2$ it is immediate that $\RXL = X$ while 
if $n\ge 3$ a generic hyperplane, tangent to $X$ at one point, 
cuts on $X$ an irreducible quadric cone and thus $\RXL = 
\emptyset.$ Let now $V$ be a spanning subspace of codimension one. 
Let us first observe that $\phiV : X \to \Pin{n}$ is a double 
cover, ramified along a smooth quadric hypersurface $\Delta = 
\mathbb{Q}^{n-1}\subset \Pin{n}.$ Let now $R$ be the ramification 
divisor of $\phiV.$ Let $x \in X.$ If $x \not \in R$ then $|V-2x|$ 
is empty, while if $x \in R$ then $|V-2x|$ consists of a unique 
divisor $D,$ corresponding to a 
 hyperplane $H,$ tangent to $\Delta$ at $\phiV(x).$ 
If $ n=2$ then $D$ is clearly reducible as the union of two lines, 
thus $x \in \RXV$ and $\RXV = R = \jonev.$ If $n \ge 
3$ then $D$ is an irreducible quadric cone and therefore $\RXL = 
\emptyset.$ Notice that $\restrict{\phiV}{D} : D \to H$ is a 
double cover ramified along the intersection of $H$ and $\Delta,$ 
which is a reducible conic. The results of this discussion are 
summarized in the following table 
 
\begin{center} 
\begin{tabular}{|c|c|c|}\hline 
$n$&${\rm codim}_{|L|}(|V|)$& $\RXL$\\ 
\hline 2 & 0 & $X$ \\ 
\hline 2 & 1 & $\jone$\\ 
\hline $\ge 3$ & 0, 1 & $\emptyset$\\ 
\hline 
\end{tabular} 
 
\end{center} 
 
\end{example} 
%
%
\begin{example}[The second symmetric product of a hyperelliptic curve] 
\label{RS4*} The following class of surfaces was considered in an 
apparently unrelated context in \cite{gisa2}. For the convenience 
of the reader we recall the construction. Let $C_1$ and $C_2$ be 
two copies of the same hyperelliptic curve of genus $q \ge 1.$ Let 
$X_q = C_1 \times C_2$ and let $\pi_i : X \to C_i,$ $i=1,2,$ be 
the projections onto the factors. Let $\iota : X_q \to X_q$ be the 
involution $\iota (x,y)=(y,x)$ and let $S$ be the quotient 
$X_q/\iota.$ Let $p:X_q \to S$ be the resulting double cover. 
 
Choose a $g^1_2$ on $C_i$ (pick the unique one if $q\ge 2$). 
Consider a divisor $D=Q_1 + Q_2\in g^1_2$ Let $\hat{L}_i=\pi^*_iD$ 
and let $\hat{L}=\hat{L}_1 + \hat{L}_2.$ It is 
$h^0(X_q,\hat{L})=4.$ Let $H^0(X_q,\hat{L})^{\iota}$ be the 
subspace of global sections of $\hat{L}$  which are 
$\iota$-invariant. If $H^0(C_i,D)=<\sigma, \tau>$ then 
$H^0(X_q,\hat{L})^{\iota}=<\sigma \tensor \sigma, \tau \tensor 
\tau, \sigma \tensor \tau + \tau \tensor \sigma>$. Now let 
$\hat{C} \in |\hat{L}|$ and consider the line bundle $L$ on $S$ 
associated to the divisor $p(\hat{C}),$  so that $p^*L=\hat{L}.$ 
There is a natural isomorphism between global sections of $L$ and 
global sections of $\hat{L}$ which are $\iota$-invariant. 
Therefore $h^0(S, L)=3.$ From the construction it follows that $L$ 
is ample and spanned, and $2L^2= \hat{L}^2 = (\hat{L}_1 + 
\hat{L}_2)^2 = 8$ so that $L^2=4.$ 
 
In \cite{gisa2}, Section 5.2, it was shown that for all these 
surfaces $\bslocus{K+L}$ contains a smooth rational component 
$\Gamma.$ Moreover, for all $x\in \Gamma$ the pencil $|L-x| $ 
contains exactly one curve $C$ singular at $x,$ which is 
reducible, while every other member of the pencil is smooth at $x$ 
and has at $x$ the same tangent direction, see \cite{gisa2} 
Section 5.2 and Proposition 6.3. The latter result already shows 
that 
\begin{equation} 
\label{GammainR} \Gamma \subseteq \RSL. 
\end{equation} 
 Now let $y_1, \dots y_{2q+2}$ be the set of the ramification points of the 
chosen $g^1_2$ on $C_i.$ These points, via the above construction, 
give rise to $(2q+2)$ reducible elements  of the form $2\Gamma_i 
\in |L|.$ 
 
Let $\mathfrak{J} = \{x_1,\dots x_s\}$, $s = \binom{2q+2}{2}$, be 
the set of the  intersection points of all pairs of $\Gamma_i's.$ 
Notice that all the $x_k$'s are distinct. If $x_k = \Gamma_r \cap 
\Gamma_s,$ a simple check in local coordinates shows that every 
element in the pencil $|L-x_k|$ is singular at $x_k$ but the only 
non reduced ones are precisely $2\Gamma_r$ and $2\Gamma_s.$ This 
also shows that for all $x_k \in \mathfrak{J},$ it is $|L - x_k| = 
|L - 2x_k|.$ Therefore it is $\mathfrak{J}\subseteq 
\mathcal{J}_2(S, L).$ Notice that the upper bound mentioned at the 
beginning of Section \ref{Bexamples} gives $\binom{2q+2}{2} = 
\text{Card}(\mathfrak{J})\le\text{Card}(\mathcal{J}_2(S,L))\le 
2q^2 + 3q + 3$ and thus $\mathfrak{J}$ misses at most two points 
of $\mathcal{J}_2(S,L).$ Notice also that 
\begin{equation} 
\label{J2notinR} 
\text{ for all  } x \in \mathcal{J}_2(S,L), \text{ 
it is } x \not \in \RSL. 
\end{equation} 
Otherwise \brref{basvjn} and the fact that $q(S) = h^1(\calo_S) > 
0$ would give $x \in \BSL$ which contradicts Proposition 
\ref{B49}. 
 
Let now $x \in \Gamma_i \setminus \mathcal{J}_2(S,L)$ for some 
$i.$ (Notice that this means $x \not \in \mathfrak{J}$  and 
possibly $x$ different from two more points.) As $x\not \in 
\mathcal{J}_2(S,L)$  and $\dim{(|L-x|)}=1,$ it is 
$\dim{(|L-2x|)}=0,$ i.e. the only singular element in $|L-x|$ is 
$2\Gamma_i$ which is reducible, and therefore 
\begin{equation} 
 \label{GammaiinR} 
\text{for all } x   \in \Gamma_i \setminus \mathcal{J}_2(S,L) 
\text { it is   } x \in \RSL. 
\end{equation} 
 
This shows that $\mathcal{W} = \Gamma \cup (\cup_{i = 
1}^{2q+2}\Gamma_i \setminus \mathcal{J}_2(S,L)) \subseteq \RSL.$ 
Recall from Lemma \ref{introtoR} that $\RSL \subseteq 
\mathcal{J}_1(S,L).$ As $h^0(L) = 3$, $\phiL$ gives a $4 -1$ cover 
of $\Pin{2}$ and thus 
 the ramification locus of $\phiL$ is a divisor, say $R.$ Let $Supp(R)$ denote the 
 reduced support of $R.$ It is $Supp(R) = \mathcal{J}_1 (S,L)$ and thus 
 
 \begin{equation} 
\label{gammasinj1} \Gamma \cup (\cup_{i=1}^{2q+2} \Gamma_i) = 
\overline{\mathcal{W}}\subseteq \overline{\RSL} \subseteq 
\mathcal{J}_1 (S,L). 
\end{equation} 
We claim that $R$ is reduced and that $Supp(R) = R = 
\mathcal{J}_1(S,L) = \Gamma + \Gamma_1 \dots \Gamma_{2q+2}$ as 
divisors on $S.$ To see this, first recall that $R$ is linearly 
equivalent to $K_{S} + 3L.$ If the claim were not true, 
\brref{gammasinj1} would imply that there exist positive integers 
$\nu, \nu_i$ and  an effective (or possibly trivial) divisor 
$\mathcal{D}$ such that 
\begin{equation} 
\label{equalityDelta} \nu \Gamma + \sum_{i=1}^{2q+2}\nu_i\Gamma_i 
+ \mathcal{D} \lineq K + 3L. 
\end{equation} 
The construction of these polarized surfaces shows that $g(L) = 
2q,$ hence $KL = 4q - 6,$ and it also shows  that $\Gamma L = 2.$ 
Computing the intersection of both sides of the last equality with 
$L$ we have 
$$2(\nu + \sum \nu_i) + L \mathcal{D} = 2 (q + 3).$$ As $\nu 
>0$ and $\nu_i > 0$ for all $i$, and $L$ is ample, the above equality implies 
$\nu = \nu_i = 1$ for all $i$ and $\mathcal{D} = \calo_X.$ 
Therefore $Supp(R) = R= \mathcal{J}_1(S,L) = \Gamma 
\cup(\cup_{i=1}^{2q+2} \Gamma_i).$ 
 
Now observe that $\restrict{\phiL}{\Gamma}$ is a $1:1$ map onto a 
smooth conic $\gamma.$ To see this, notice that the only other 
possibility would be for $\restrict{\phiL}{\Gamma}$ to be a $2:1$ 
cover of a line. In this case, let $z$ be a general point on 
$\gamma.$ Then $\restrict{\phiL}{\Gamma}^{-1}(z) = \{x,x'\}$ where 
$x,x' \in \Gamma, x \neq x'.$ It follows that $|L - x| = | L - 
x'|$ but this contradicts the accurate description of $|L-x|$ 
given in \cite{gisa2}, Proposition 6.3. In particular, the only 
singular element $D \in |L-x|$ is of the form $D=A+B$ with $A 
\Gamma = B \Gamma = 1$ and thus $ D \not \in |L-x'|.$ 
 
Further observe that $\restrict{\phiL}{\Gamma_i}$ is a $2:1$ map 
onto a line $\gamma_i,$  tangent to $\gamma,$ for all $i=1, 
\ldots, 2q+2.$ To see this it is enough to recall that $L= 
2\Gamma_i$ and that $\Gamma \Gamma_i = 1$ for all $i.$ 
 
The above discussion on the images of the components of the 
ramification divisor of $\phiL$ shows that the branch locus of 
$\phiL$ consists of the union of $\gamma$ with $2q+2$ tangent 
lines $\gamma_i.$ 
 
Then Lemma \ref{J2andSing} and the fact that $\cup_i ( \Gamma \cap 
\Gamma_i) \not \subset \mathcal{J}_2(S,L)$ gives 
\begin{equation} 
\label{J2ofS4*} \mathcal{J}_2(S,L)= \{x_1 ,\ldots, x_{2q+2}\}. 
\end{equation} 
 
Finally, \brref{GammainR}, \brref{J2notinR}, \brref{GammaiinR} and 
\brref{J2ofS4*} give: 
 
$$\RSL = 
\Gamma \cup (\cup_{i=1}^{2q+2} \Gamma_i \setminus \{x_1 ,\ldots, 
x_{2q+2}\}).$$ 
\end{example} 
%
%
 
\begin{example}[Double cover of $\Pin{1} \times \Pin{1}$] 
\label{RdoublecoverofQ2} Let $a \geq 2$ be an integer and let $S$ 
be the smooth surface defined by the double cover $\pi:S \to 
\mathbb{Q}$ of the smooth quadric surface $\mathbb{Q} =\Pin{1} 
\times \Pin{1}$ branched along a smooth curve $\Delta \in 
|\mathcal O_{\mathbb{Q}}(2a,2a)|$. Let $L:= \pi^*\mathcal 
O_{\mathbb{Q}}(1,1)$. Then $L$ is an ample and spanned line bundle 
on $S$. By using the projection formula we see that 
$$h^0(L)=h^0(\pi_*L)=h^0(\mathcal{O}_{\mathbb{Q}}(1,1)+ 
h^0(\mathcal O_{\mathbb{Q}}(1-a,1-a))=4.$$ This shows that the 
morphism $\varphi_L:S \to \Pin{3}$ factors through $\pi$; in 
particular, $L$ is not very ample. Let $x \not\in \mathcal 
{J}_1(S,L)$. In this case $|L-2x| = \{\pi^*h\}$, where $h$ is the 
only element in $|\mathcal{O}_{\mathbb{Q}}(1,1)-2\pi(x)|$. Thus 
$D$ is reducible, so being $h$, the section cut out on 
$\mathbb{Q}$ by its tangent plane at $\pi(x)$. Therefore $x \in 
\mathcal{R}(S,L)$. Now let $x \in \mathcal{J}_1(S,L)$. Recall that 
$\mathcal{J}_1(S,L)$ is the ramification curve of the double cover 
$\pi$, hence it is isomorphic to the branch curve $\Delta$, via 
$\pi$. Thus 
$$|L-2x|=\{D = \pi^*h \quad | \quad h 
\in |\mathcal{O}_{\mathbb{Q}}(1,1)-\tau_{\pi(x)}|\},$$ where 
$\tau_{\pi(x)}$ is a tangent vector to $\Delta$ at $\pi(x)$.  Note 
that the general element $h$ as above is irreducible, since it is 
cut out on $\mathbb{Q}$ by a plane of $\Pin{3}$ not tangent to 
$\mathbb{Q}$ itself. Hence the corresponding element $D$ in the 
pencil $|L-2x|$ is also irreducible. This means that $x \not\in 
\mathcal{R}(S,L)$. In conclusion we have: 
$$\mathcal{R}(S,L) = S \setminus \mathcal{J}_1(S,L).$$ 
\end{example} 
 
\section{Towards the Conjecture} 
 
The examples in Section \ref{Rexamples} lead us to formulate the 
following conjecture. 
 
\begin{conj} 
\label{congettura} 
 Let $X$ be a smooth complex variety of 
dimension $n\ge 2.$ Let $L$ be an ample line bundle on $X,$ 
spanned by a subspace $V \subseteq H^0(L)$ with $\dim{(|V|)} \ge 
n+1.$ If $\RXV = X$ then $\xel$ is either $(\Pin{2}, 
\oofp{2}{2})$, $(\Pin{1} \times \Pin{1}, \oof{\Pin{1} \times 
\Pin{1}}{1,1})$ or a scroll over a smooth curve. 
\end{conj} 
 
Conjecture \ref{congettura} is very easily verified if one 
strengthens the requirements on $|V|$ by asking for its very 
ampleness. One can even relax a bit the condition on the size of 
$\RXV.$ 
 
\begin{prop} 
Let $X$ be a smooth complex variety of dimension $n\ge 2.$ Let $L$ 
be a line bundle on $X,$ with $V \subseteq H^0(L)$ such that $|V|$ 
is very ample and $\dim{(|V|)} \ge n+1.$ If $\RXV$ is Zariski 
dense in $X,$ then $\xel$ is either $(\Pin{2}, \oofp{2}{2})$, 
$(\Pin{1} \times \Pin{1}, \oof{\Pin{1} \times \Pin{1}}{1,1})$ or a 
scroll over a smooth curve. 
\end{prop} 
\begin{proof} 
Consider $X$ as embedded by $|V|.$ Let $x$ be a general point in 
$X$ and thus $x \in \RXV.$ If $n=2$ the general hyperplane tangent 
to $X$ at $x$ has no other tangency locus with $X$ and thus it 
cuts on $X$ a hyperplane section $D \in |L-2x|$ that is reducible 
and has a single ordinary quadratic singularity at $x.$ The proof 
now proceeds exactly as \cite{BESO}, Corollary 1.6.8, p.\ 31. 
 
Let $n\ge 3.$  For all $D \in |L-2x|$, $D$ is reducible and thus 
\linebreak $\dim{(\text{Sing}(D))}~\ge~n-2.$ Notice that 
$\text{Sing}(D)$ is the contact locus of the hyperplane 
corresponding to $D$ and $X.$ This being true for a general $x \in 
X$ implies that the dual defect of $(X,V)$ is $\geq n-2.$ The 
conclusion now follows from \cite{la-stru}, Corollary 3.4 or 
\cite{Ein1}, Theorem 3.2. \qed \end{proof} 
 
\begin{rem}Example \ref{RdoublecoverofQ2} shows that the very ampleness is 
necessary if the size of $\RXV$ is relaxed, as there $\RXL$ is 
Zariski dense, $L$ is ample and spanned but not very ample. 
\end{rem}  
%
%

\end{document}